\documentclass[11pt]{amsart}

\newif\ifdraft
\draftfalse

\usepackage{amsmath,amsthm,verbatim,amssymb,amsfonts,amscd, graphicx}
\usepackage{graphics}
\usepackage{theoremref}
\usepackage{mathtools}
\usepackage[alphabetic,initials]{amsrefs}

\usepackage{tikz-cd}

\usepackage{amsmath, amscd}%
\usepackage{amsfonts}%
\usepackage{amssymb}

\usepackage{graphicx}
\usepackage[pdftex=true, plainpages=false]{hyperref}

\renewcommand{\tilde}{\widetilde}
\newtheorem{thm}{Theorem}[section]

\newtheorem{cor}[thm]{Corollary}

\newtheorem{lem}[thm]{Lemma}
\newtheorem{prop}[thm]{Proposition}

\theoremstyle{remark}
\newtheorem{ex}[thm]{Example}

\theoremstyle{definition}
\newtheorem{de}[thm]{Definition}
\newcommand{\Dmod}{\mathfrak{D}}

\newcommand{\LL}{{\mathcal{L}}}

\newcommand{\R}{{\mathbb R}}

\newcommand{\C}{{\mathbb C}}
\newcommand{\Q}{{\mathbb Q}}
\newcommand{\Z}{{\mathbb Z}}

\newcommand{\OO}{{\mathcal{O}}}

\newcommand{\map}{\rightarrow}

\newcommand{\xmap}{\xrightarrow}

\DeclareMathOperator{\End}{End}
\DeclareMathOperator{\Ker}{Ker}
\DeclareMathOperator{\var}{var}
\DeclareMathOperator{\Var}{Var}

\DeclareMathOperator{\can}{can}

\DeclareMathOperator{\Tr}{Tr}
\DeclareMathOperator{\im}{Im}

\DeclareMathOperator{\Cone}{Cone}
\DeclareMathOperator{\gr}{Gr}
\DeclareMathOperator{\DR}{DR}

\newcommand{\be}{\begin{enumerate}}
\newcommand{\ee}{\end{enumerate}}
\newcommand{\bt}{\begin{thm}}
\newcommand{\et}{\end{thm}}
\newcommand{\bde}{\begin{de}}

\newcommand{\ede}{\end{de}}
\newcommand{\bc}{\begin{cor}}
\newcommand{\ec}{\end{cor}}
\newcommand{\blm}{\begin{lem}}
\newcommand{\elm}{\end{lem}}
\newcommand{\bp}{\begin{proof}}
\newcommand{\ep}{\end{proof}}
\newcommand{\beq}{\begin{equation*}\label{xx}}
\newcommand{\eeq}{\end{equation*}}

\renewcommand{\nu}{{\mathcal{V}}}

\newcommand{\hm}{\mathcal{H}om}

\newcommand{\arr}{\longrightarrow}

\renewcommand{\ref}[1]{\hyperref[#1]{\ref*{#1}}}

\renewcommand{\d}[1]{{#1^\bullet}}

\newcommand{\nearby}[1]{\Psi_f\d{#1}}
\newcommand{\nearbyd}[2]{\Psi_{f, #2}\d{#1}}
\newcommand{\vanish}[1]{\Phi_f\d{#1}}
\newcommand{\vanishd}[2]{\Phi_{f, #2}\d{#1}}
\newcommand{\phinity}{\phi_{\infty}}

\newcommand{\biG}[1]{G(#1, \d K)}

\begin{document}
\title{Nearby and Vanishing cycles for Perverse Sheaves and D-modules}

\author{Lei Wu}
\address{Department of Mathematics, University of Utah,
155 S 1400 E, Salt Lake City, UT 84112, USA}
\email{\tt lwu@math.utah.edu}
\date{}

\begin{abstract}
We survey nearby and vanishing cycles for both perverse sheaves and $\Dmod$-modules under analytic setting. Following ideas of A. Beilinson, M. Kashiwara and M. Saito, we explain in detail the proof of the comparison theorem between them in the sense of Riemann-Hilbert correspondence.  
\end{abstract}
\maketitle
\section{Introduction}
The functors of nearby and vanishing cycles was first introduced by A. Grothendieck (see \cite{SGA7}). They are widely studied from different points of views, for instance by Beilinson \cite{BGluePer} algebraically, and Kashiwara and Schapira \cite{KSSheavesM} from microlocal perspectives. Gabber showed that they preserve perversity (up to a shift of degrees); see Theorem \ref{nvperverse} below. 

Passing to the category of $\Dmod$-modules, as suggested by Riemann-Hilbert correspondance (\cite[Theorem 5.7]{Kash}), there is a notion of nearby and vanishing cycles for regular holonomic $\Dmod$-modules due to Kashiwara and Malgrange involving the use of what is now called the Kashiwara-Malgrange filtration. The comparison theorem between nearby and vanishing cycles for perverse sheaves and regular holonomic $\Dmod$-modules is established in \cite[Theorem 2]{KashVan}. A more refined version of the comparison is proved in \cite[\S 3.4]{SaitoHMP} but for special $\Q$-specializable holonomic $\Dmod$-modules (what the author called holonomic $\Dmod$-modules that are quasi-unipotent and regular along a smooth hypersurface). The core  of this article is to survey the proof of the comparison theorem in $loc.$ $cit.$. We use local systems given by the infinite dimensional Jordan blocks and their Deligne canonical extensions to understand nearby cycles for both perverse sheaves and $\Dmod$-modules, which help us simplify the extremely complicated arguments used in $loc.$ $cit.$.

In \S 2, we discuss general properties of Kashiwara-Malgrange filtrations, essentially due to Kashiwara. In \S 3, we recall the construction of nearby and vanishing cycles in general following \cite[\S 8.6]{KSSheavesM}. We also construct $\lambda$-nearby cycles for perverse sheaves alternatively by using  the infinite dimensional Jordan blocks inspired by ideas in \cite{BGluePer} in this section; that is the content of Theorem \ref{altnearby}. \S 4 is about definitions of nearby and vanishing cycles for specializable $\Dmod$-modules along arbitrary hypersurfaces and the proof of the comparison theorem (see Theorem \ref{fcompnv}). Parallel to the perverse case, we also give a description of nearby cycles for specializable $\Dmod$-modules via Deligne canonical extensions given by infinite dimensional Jordan blocks, on which the proof of the comparison strongly relies; see Corollary \ref{limcomp}. 

\noindent 
{\bf Acknowledgement.} 
I thank my advisor Mihnea Popa for holding a Hodge Module learning seminar at Northwestern University and giving me a chance to lecture and learn the beautiful topic about nearby and vanishing cycles. I am also grateful to Yajnaseni Dutta. The first draft of this article is based on her notes from my lectures. I also thank Takahiro Saito for pointing out a mistake in an early version of this paper.

\section{Kashiwara-Malgrange Filtrations on $\Dmod$-modules} In this section, we introduce the Kashiwara-Malgrange filtrations and prove the existence of Kashiwara-Malgrange filtrations for specializable $\Dmod$-module. All $\Dmod$-modules are assumed to the left ones in this article. 

Let $X$ be a complex manifold of dim $n$, $H\subset X$ a smooth hypersurface and $I_H$ ideal sheaf of $H$.  
\begin{de}[Kashiwara-Malgrange filtration]
The \textit{Kashiwara-Malgrange} filtration on $\Dmod_X$ is defined by
$$V^i\Dmod_X= \{P\in \Dmod_X| PI^j_H\subset I^{j+i}_H \text{for all }j\in\Z\}$$
where $I^j_H = \OO_X$ for $j\leq 0$. 
\end{de}
Locally on coordinates $(z_1,\dots,z_{n-1},t)$ of $X$, if $H$ is defined by $t=0$, then we have
$$V^0\Dmod_X= \OO_X<\partial_1,\dots,\partial_{n-1},t\partial_t>,$$ 
$V^i\Dmod= t^iV^0\Dmod_X$ and $V^{-i}\Dmod_X= \partial_tV^{-i+1}\Dmod_X+V^{-i+1}\Dmod_X$ for $i>0$.

\noindent
\begin{de}
The $Kashiwara$-$Malgrange$ $filtration$ along $H$ on a left $\Dmod_X$-module $M$ is a $\Z$-indexed decreasing filtration $\d{V}M$ such that
\be
\item $\bigcup V^kM=M$ and $V^kM$ is coherent over $V^0\Dmod_X$;
\item $tV^kM\subset  V^{k+1}M$, $\partial_t^kM\subset V^{k-1}M$;
\item $tV^kM = V^{k+1}M$ for $k\gg 0$;
\item the $t\partial _t$-action on $\gr_{V}^k{M}=\frac{V^kM}{V^{k+1}M}$ has a minimal polynomial locally (globally in algebraic case), and eigenvalues of $t\partial_t$ have real parts in $[k, k+1)$.
\ee
If a filtration $\Omega^\bullet M$ only satisfies Condition (1), (2) and (3), then we call it a coherent filtration with respect to $(\Dmod_X,V^\bullet)$. For coherent $\Dmod_X$-modules, Kashiwara-Malgrange filtrations may not exist; but coherent filtrations always exist at least locally; see for instance \cite[Appendix A.1]{Kash}.
\end{de}
\blm
The Kashiwara-Malgrange filtration is unique if it exists.
\elm
\bp
Suppose $V^\bullet M$ and $U^\bullet M$ are two filtrations satisfying all the conditions. By symmetry, it is enough to prove that $U^k\subset V^k$
for all $k$. By coherence, locally on a neighborhood 
\[U^kM\subset V^jM\]
for $j\ll0$. 

From condition (4) we know that there exists $b_1(s)\in \C[s]$ whose roots have real parts in $[k, k+N)$ such that $b_1(t\partial_t)$ annihilates $\frac{U^kM}{U^kM\cap V^{j+1}M}$ by using the surjection
\[\frac{U^kM}{U^{k+N}M}\twoheadrightarrow \frac{U^kM}{U^kM\cap V^{j+1}M}\]
$N\gg0$. On the other hand, there exists $b_2(s)\in \C[s]$ whose roots have real parts in $[j, j+1)$ such that $b_2(t\partial_t)$ annihilates $\frac{U^kM}{U^kM\cap V^{j+1}M}$ by using the injection
\[\frac{U^kM}{U^kM\cap V^{j+1}M}\hookrightarrow \frac{V^jM}{V^{j+1}M}.\] If $j<k$, since $b_1(s)$ and $b_2(s)$ have no common root, then by Bez\'out's lemma 1 kills $\frac{U^kM}{U^kM\cap V^{j+1}M}$. Hence $U^kM\subset V^{j+1}M$.
Repeating this process, we obtain $U^kM\subset V^{k}M$. 
\ep

\begin{de}[$Specializability$]\label{smspe}
 A coherent $\Dmod_X$-module $M$ is said to be specializable along $H$ if locally for some coherent filtration $\d{\Gamma}$ there exists a $monic$ polynomial $b(s)=b_\Gamma(s)\in\C[s]$ such that $b(t\partial_t - k)$ acts on $\gr_k^\Gamma M$ trivially. 
 
If we use $R$ to denote the field $\Q$ or $\R$, then $M$ is $R$-$specializable$ if additionally roots of $b(s)$ are contained in $R$. For instance, Hodge modules are $\Q$-$specializable$; see \cite{SaitoHMP}.
\end{de}

\begin{lem}
Suppose $M$ is specializable along $H$. Then locally there exists a coherent filtration $\d{\Omega}M$ of $M$ such that there is a polynomial $b'(s)\in\C[s]$ satisfying the condition that $b'(t\partial_t - k)$ acts on $\Omega^kM/\Omega^{k+1}M$ trivially and the real parts of  $\text{roots of } b'(s)\subset [0,1)$.
\end{lem}
\bp
Set $\Omega^n = \Gamma^{n+k}$. Then $b(t\partial_t - n-k)$ acts on $\Omega^n/\Omega^{n+1}$ trivially. So we can assume that the real part of the roots of $b(s)$ for $\Gamma^\bullet$ are larger than 0. 

Assume $b(s)=(s-\alpha_1)^{n_1}b_1(s)$. Now set $\Omega^n = \Gamma^{n+1} + (t\partial_t - \alpha_1-n)^{{n_1}}\Gamma^n$. We can see that $(t\partial_t - (\alpha-1)-n)^{n_1}b_1(t\partial_t-n)$ annihilates $\Omega^n/\Omega^{n+1}$. Repeating this process, we can move the real parts of roots of $b(s)$ all in $[0,1)$.
\ep
By uniqueness of Kashiwara-Malgrange filtrations, the above lemma immediately implies the following corollary.
\bc\label{speKM}
If $M$ is specializable along $H$, then the Kashiwara-Malgrange filtration of $M$ along $H$ exists globally. 
\ec

The definition of $specializability$ is due to Sabbah. From the above corollary, it is equivalent to the existence of Kashiwara-Malgrange filtrations. Indeed, they are also equivalent to the existence of (generalized) $b$-functions; see \cite{SabVan}.

\begin{ex}[Kashiwara's equivalence]
If $M$ is supported on $H$, then $M$ is specializable along $H$. In this case, the existence of Kashiwara-Malgrange filtration is equivalent to Kashiwara's equivalence. Moreover,
\[M\simeq i_+\gr^{-1}_VM,\]
where $i: H\hookrightarrow X$ is the closed embedding and $i_+$ is the $\Dmod$-pushforward. See \cite[Theorm 1.6.1]{HTT} or \cite{Kash} for details.
\end{ex}

\bt[M. Kashiwara]
If $M$ is holonomic, then $M$ is specialisable along every smooth hypersurface $H\subset X$. Consequently, its Kashiwara-Malgrange filtration exists and $\gr^k_V M$ is holonomic. Moreover, if $M$ is also regular, then so is $\gr^k_V M$ on $H$.
\et
Let us refer to \cite{KashVan} for the definition of (regular) holonomic $\Dmod$-modules and a proof of the above fundamental theorem; see also \cite{Bj} for a more algebraic approach.

\section{Nearby and Vanishing Cycles for Perverse Sheaves} We will discuss nearby and vanishing cycles for perverse sheaves in this section following \cite[\S 8.6]{KSSheavesM}. Through out this section,  $X$ will be a complex manifold of $\dim n$. 

\textit{\bf Notations.}  If $f$ is a morphism of complex manifolds, we use notations $Rf_*$ the derived push-forward, $f^{-1}$ the sheaf pullback, $Rf_!$ the derived push forward with compact support, $f^!$ the adjoint functor of $Rf_!$.. 

\subsection{Decompositions of nearby cycles and vansihing cycles} 
For a holomorphic function $f$ on $X$, consider the following morphisms:
\[\begin{tikzcd}
&&\tilde{X}\arrow{r}\arrow{d}{p}&\tilde{\C}^*\arrow[d, "\pi"]\\
X_0\arrow[r, hook, "i"]&\arrow[ur, leftarrow, "\tilde{p}"] X\arrow[r, hookleftarrow, "j"]  \arrow[rr, bend right, "f"]& X^*  \arrow[r, "f_0"]&\C
\end{tikzcd}\]
where $X^*  = f^{-1}(\C^*)$ and $X_0=f^{-1}(0)$, and where $\pi$ is the composition of the universal covering map of $\C^*$ and the open embedding $j:\C^*\hookrightarrow \C$. For a constructible complex $\d{K}$ on $X$, we recall the definition of nearby cycle, 
$$\nearby{K} = i^{-1}Rj_*Rp_*p^{-1}j^{-1}\d{K}=i^{-1}R\tilde{p}_*\tilde{p}^{-1}\d{K}.$$ From definition, it only depends on $j^{-1}\d{K}$.

By Poincar\'e-Verdier duality, we know
$$\nearby{K}=i^{-1}R\hm_{{\C_{X}}}(f^{-1}{\pi}_!\C_{\tilde{\C^*}}, \d{K}).$$ Since the fibre of $\pi$ is a discrete set and isomorphic to $\Z$, $\pi_!\C_{\tilde{\C^*}}|_{\C_*}$ is a $\C[\Z]$-local system of rank 1. Let $T$ be the operator of the monodromy action on $\pi_!\C_{\tilde{\C^*}}$ around the origin counterclockwise. The monodromy action induces an action on the nearby cycle, also denoted by $T$.
Consider diagram $\Tr_{\can}$,
\begin{equation}\label{Scan}
\Tr_{\can}:
\begin{CD}
0@>>>0@>>>\pi_!\C_{\tilde{\C^*}} @>>id>\pi_!\C_{\tilde{\C^*}} @>>>0\\
@.@VVV @VtrVV @VVV\\
0@>>>\C_{\C}@>>>\C_{{\C}} @>>>0 @>>> 0
\end{CD}
\end{equation}

The trace map is defined via adjunction. Since
\[
    (\pi_!\C_{\tilde{\C^*}})_a=\left\{   
                         \begin{array}{ll}
                               0, \textup{ if } a=0\\
                               \C[\Z]\simeq\bigoplus_{\Z}\C, \textup{ otherwise},                
                                  \end{array}
              \right.
  \]
stalkwise the trace map is the same as taking sum of all the entries. Treat the vertical maps as complexes $\d{A}$, $\d{B}$ and $\d{C}$ respectively in which $\pi_!\C_{{\C^*}}$ is of degree 0.  It is clear that $\Tr_{\can}$ is a short exact sequence of complexes.
Set $G(\bullet,\bullet)=i^{-1}R\hm_{\C_{{X}}}(f^{-1}\bullet,\bullet)$. The vanishing cycle is defined by 
$$\vanish{K}=\biG{\d{B}}.$$
Appying $\biG{\bullet}$ to $\Tr_{\can}$ we get,
\begin{equation}\label{canses}
 i^{-1}\d{K}\map\nearby{K} \xmap{\can} \vanish{K} 	 \xmap{+1}.
\end{equation}
Namely $\vanish{K}$ is the cone of the natural map $i^{-1}\d{K}\map\nearby{K}$.
Now consider another short exact sequence of complexes $\Tr_{T-1}$
\begin{equation}
\Tr_{T-1}:\begin{CD}
0@>>>\pi_!\C_{\tilde{\C^*}}@>T-I>>\pi_!\C_{\tilde{\C^*}} @>tr>>j'_{!}\C_{{\C^*}} @>>>0\\
@.@VVV @VtrVV @VVV\\
0@>>>0@>>>\C_{{\C}} @>id>>\C_{{\C}} @>>> 0
\end{CD}\end{equation}
Applying $\biG\bullet$, we get another triangle 
\begin{equation}\label{varses}
i^!\d{K}\map\vanish{K} \xmap{\var} \nearby{K}\xmap{+1}.
\end{equation}

We have constructed $\can$ and $\var$ between the nearby cycle and the vanishing cycle. The monodromy action on $\d B$ (T acts trivially on $\C_C$) induces an $T$-action on  $\vanish{K}$ . By construction, $\can\circ\var= T-I$ and $\var\circ\can = T-I$. 

Now we look at the eigenvalue decomposition of nearby and vanishing cycles with respect to the $T$-action. First, we need the following result due to O. Gabber. See for instance \cite{Bry} for more information.
\bt \label{nvperverse}
If $\d K$ is perverse, then $\nearby{K}[-1]$ and $\vanish{K}[-1]$ are also perverse. 
\et

From now on, $\d K$ is assumed to be perverse. We know the category of perverse sheaves is abelian (\cite[\S 8]{HTT}). The above theorem implies that triangles $\biG{\Tr_{\can}}$ and $\biG{\Tr_{1-T}}$ give rise to two exact sequences in the abelian category after taking perverse cohomologies. 

The following lemma working abstractly for abelian categories is also needed.

\blm\label{decomp}
Suppose $A$ is an object in a $\C$-abelian category ($\textup{Hom}$'s are $\C$-vector spaces) with an isomorphism $\varphi$. If $g(\varphi)=0$ for some $g(x)\in \C[x]$, then $A$ has a unique generalized eigenspace decomposition (the Jordan decomposition) with respect to the $\varphi$-action. 
\elm

\bp
Suppose $\lambda$ is a root of $g(x)$ with multiplicity $m$ and $A_\lambda=\ker(\varphi-\lambda)^m$. Then we have a short exact sequence
\[0\map A_\lambda\map A\map Q\map 0,\]
where $Q=\im{(\varphi-\lambda)^m}$ . However $(\varphi-\lambda)^m: Q\map A$ splits the short exact sequence which proves the statement. 

The decompostion is unique because each factor is universally defined.
\ep

Clearly, if $m'>m$, $\ker(\varphi-\lambda)^{m'}=\ker(\varphi-\lambda)^m$. Hence, it makes sense have the identity $A_\lambda= \ker(\varphi-\lambda)^\infty$.
\begin{prop}
If $\d K$ is perverse, then $\nearby{K}$ and $\vanish{K}$ have functorial decompositions with respect to the $T$-action,
\[\nearby{K}\simeq \bigoplus_{\lambda\in \C^*}\Psi_{f, \lambda}\d{K},\]
and 
\[\vanish{K}\simeq \bigoplus_{\lambda\in \C^*}\Phi_{f, \lambda}\d{K}.\] 
\end{prop}
\bp
We only prove the statements for vanishing cycles. That for nearby cycles can be proved similarly. 

Locally on a relative compact open neighborhood, the $T$-action on $\vanish{K}$ clearly has a polynomial $g(x)\in \C[x]$ such that $g(T)=0$. Since the category of perverse sheaves is an abelian category, by the above lemma, $\vanish K$ has a decomposition at least locally. On the other hand, globally $\Phi_{f,\lambda}{\d K}=\ker{(T-\lambda)^\infty}$. Therefore, the decomposition is global.
 
Assume $\varphi$ is a morphism of perverse sheaves on $X$. We know the $T$-action on $\Phi_{f}(\bullet)$ is induced from the $T$-action on $ f^{-1}([\pi_!\C_{\tilde{C^*}}\map \C_\C])$ which stands on the first entries of $R\hm_{\C_{{X}}}(\bullet, \bullet)$. Therefore, $\Phi_f(\varphi)$ and $T$ commute, from which functoriality follows. 
\ep

From the proof, we see the decompositions are locally finite. Also when $T$ is locally quasi-unipotent (that is $(T^m-1)^n=0$), $\lambda$ can only be roots of unity; in particular if $T$ is unipotent, then $\Psi_{f, 1}\d K=\nearby K$ and $\Phi_{f, 1}\d{K}=\vanish K$.  

The $\lambda$-nearby and $\lambda$-vanishing cycles possess much more abundant structures. In fact, when $\d K$ underlies a polarizable Hodge modules, $\Psi_{f, \lambda}{\d K}$ and $\Phi_{f, 1}\d{K}$ underly mixed Hodge modules \cite{SaitoMHM}.

With the help of these decompostions, we can refine triangles \eqref{canses} and \eqref{varses}. To be precise, by construction $T$ acts on triangles $\biG{\Tr_{\can}}$ and $\biG{\Tr_{1-T}}$; hence they also have the generalized eigenspace decompositions. In particular, $\can$ and $\var$ decompose accordingly.  In summary, we obtain the following theorem. 
\begin{thm}\label{unipvarses}
If $\d K$ is perverse, then we have the followings.
\begin{enumerate}

\item For each $\lambda$, $\can$ and $\var$ induce morphisms
\[\can: \Psi_{f,\lambda}\d K\map \Phi_{f,\lambda}\d K\]
and 
\[\var: \Phi_{f,\lambda}\d K\map \Psi_{f,\lambda}\d K.\]
\item $\var\circ\can=T-I$ and $\can\circ\var=T-I$. 
\item If $\lambda\neq 1$, $\can: \Psi_{f,\lambda}\d K\map \Phi_{f,\lambda}\d K$ and $\var: \Phi_{f,\lambda}\d K\map \Psi_{f,\lambda}\d K$ are isomorphisms.
\item There are triangles
\[ i^{-1}\d{K}\map \nearbyd{K}{1} \xmap{\can} \vanishd{K}{1} \xmap{+1};\]
and
\[i^!\d{K}\map\Phi_{f, 1} \d{K} \xmap{\var} \Psi_{f,1}\d{K}\xmap{+1}.\]
\end{enumerate}
\end{thm}

Since the canonical morphism $\can$ is isomorphic, there is no need to make a distinction between $\Psi_{f,\lambda}\d K$ and $\Phi_{f,\lambda}\d K$ for $\lambda\neq 1$; that is we indentify
\begin{equation} \label{altlambda}
\Psi_{f,\lambda}\d K\stackrel{\can}{\simeq}\Phi_{f,\lambda}\d K
\end{equation}
for $\lambda\neq 1$ and hence also
\begin{equation}\label{altone}
\vanishd{K}{1}\stackrel{\can}{\simeq} \Cone(i^{-1}\d K\map \Psi_{f,1}\d K).
\end{equation} 

The operator $T$ on $\nearby{F}$ and $\vanish{F}$ have the Jordan-Chevalley decomposition
\[T=T_s\circ T_u=T_u\circ T_s,\]
where $T_s$ is semi-simple, and $T_u$ is unipotent. Indeed, the decompositions exist locally on relative compact neighborhood and local decompositions glue by uniqueness. 

Besides the morphism $\can$ and $\var$, there is another morphism "$\Var$" originally constructed in \cite{KashVan}. Let us give the definition in our setting. Under the above identifications \eqref{altlambda} and \eqref{altone}, the morphism $\Var_\lambda: \Phi_{f,\lambda}\d K\arr\Psi_{f,\lambda}\d K$ is defined by 
\[\Var_\lambda=
\begin{cases}
(0, \log T_u) &\text{ if } \lambda=1,\\
\log T_u &\text{ otherwise.}
\end{cases}\]
Then the total $\Var$ is
\[\Var:=\bigoplus_\lambda \Var_\lambda :\vanish{ K}\map \nearby{ K}.\]
Immediately from definition, we have the following:
\begin{prop}
We have $\Var\circ\can=\log T_u$ and $\can\circ\Var=\log T_u$.
\end{prop}

\subsection{Nearby cycles via local systems on $\C^*$}
In this section, we will give an alternative approach to understand nearby and vanishing cycles via local systems on $\C^*$.

\subsubsection{Local systems on $\C^*$}
Denote the counterclockwise loop around the origin in $\pi_1(\C^*)$ by $T$.  Suppose $L$ is a regular local system over $\C$ (of finite type) on $\C^*$; that is, the matrix of the monodromy $T$ is regular. By classical Riemann-Hilbert correspondence, $L$ is identified with the pair $(L_x, T)$ where $T$ is the monodromy action on $L_x$ for $x\in \C^*$. 

For any $g(T)\in C[T]$, define a local system $L_g$ by a pair $(\dfrac{\C[T]}{(g(T))}, T)$ where the $T$-action is multiplication by $T$.

\begin{lem}
Suppose $L$ is a regular local system and $g$ is the characteristic polynomial of the monodromy action $T$. Then 
\[L\simeq L_g\]
\end{lem}
\bp
Since the isomorphism class of a regular local system on $\C^*$ is uniquely determined by its characteristic polynomial, we know $L\simeq L_g$.
\ep

\begin{ex}[Local systems given by Jordan blocks]\label{lsjb}
For $\alpha\in \C$ and $\lambda=e^{2\pi\sqrt{-1}\alpha}$, we have 
\begin{equation}\label{JordanLS}
L_{(T-\lambda)^m}\simeq H_{\alpha, m}
\end{equation} 
where $H_{\alpha, m}$ is the rank $m$ local system with the monodromy action by matrix $e^{2\pi\sqrt{-1}J_{\alpha,m}}$ and $J_{\alpha, m}$ the $m\times m$ Jordan block with eigenvalues equal to $\alpha$. 
\end{ex}

\blm \label{dualrr}
We have 
$$L_{(T-\lambda)^m}^\vee\simeq L_{(T-\lambda^{-1})^m}.$$
The monodromy action on $L_{(T-\lambda)^m}$ induces an action on its dual $L_{(T-\lambda)^m}^\vee$ which is the inverse of the monodromy action on $L_{(T-\lambda^{-1})^m}$ under this isomorphism.
\elm
\bp
Both statements are obvious, because $L_{(T-\lambda)^m}^\vee$ is regular.
\ep
\blm
For each $g(T)\in \C[T]$ satisfying $g(0)\neq 0$, there exists a short exact sequence $S_g$,
\begin{equation}\label{mT}
\begin{CD}
S_g: 0@>>>\pi_!\C_{\tilde{\C^*}}@>g(T)>>\pi_!\C_{\tilde{\C^*}} @>tr_g>>j'_!L_g @>>>0.
\end{CD}
\end{equation}
Moreover, the monodromy $T$ acts on $S_g$.
\elm
\bp
Clearly, the local system $\pi_!\C_{\tilde{\C^*}}|_{\C^*}$ is given by $\C[T,T^{-1}]$ with $T$-action via multiplication by $T$. Therefore, the restriction of $S_g$ on $\C^*$ can be represented by 
\[\begin{CD}
0@>>>\C[T,T^{-1}]@>g(T)>>\C[T,T^{-1}] @>>>\frac{\C[T]}{(g(T))}@>>>0.
\end{CD}\]
The natural projection $\C[T,T^{-1}] \map \frac{\C[T]}{(g(T))}$ induces the morphism $tr_g$.
Then both of the assertions are clear.
\ep
When $g(T)=T-I$, the $S_{T-I}$ recovers
\[\begin{CD}
0@>>>\pi_!\C_{\tilde{\C^*}}@>T-I>>\pi_!\C_{\tilde{\C^*}} @>tr>>j'_!\C_{{\C^*}} @>>>0.
\end{CD}\]

\begin{thm}\label{altnearby}
Suppose $\d K$ is perverse.The morphism $tr_{(T-\lambda)^m}$ for $m\gg0$ induces an isomorphism
\[tr^{\infty}_\lambda : \varinjlim_m i^{-1}Rj_*(j^{-1}\d K\otimes f_0^{-1}H_{-\alpha,m})\map \Psi_{f,\lambda}(\d K),\]
where $\lambda=e^{2\pi\sqrt{-1}\alpha}$.  Under this isomorphism, the monodromy action $T$ on the $\lambda$-nearby cycle is induced by the action of matrix $e^{-2\pi\sqrt{-1}J_{-\alpha,m}}$ on $H_{-\alpha,m}$.
\end{thm}
\bp
For $m>0$, the $T$-action on $S_{(T-\lambda)^m}$ induces $T$-action on $\biG{S_{(T-\lambda)^m}}$. Since the $T$-action on $\biG{L_{(T-\lambda)^m}}=i^{-1}Rj_*(L^\vee_{(T-\lambda)^m}\otimes \d K)$ has only 1 eigenvalue $\lambda$, we have a distinguished triangle
\[
D_m:\begin{CD} i^{-1}Rj_*(j^{-1}\d K\otimes f_0^{-1}L^\vee_{(T-\lambda)^m})@>{tr_{(T-\lambda)^m}}>> \nearbyd{K}{\lambda} @>{(T-\lambda )^m} >> \nearbyd{K}{\lambda} @>{+1}>>,
\end{CD}
\]
which is the (generalized) $\lambda$-eigenspace of $\biG{S_{(T-\lambda)^m}}$. 

The multiplication by $T-\lambda$ induces a natural morphism of complexes
\[T-\lambda: S_{(T-\lambda)^{m+1}}\arr S_{(T-\lambda)^{m}};\]
to be precise, the morphism is
\[
\begin{CD}
0@>>>\C[T,T^{-1}]@>(T-\lambda)^{m+1}>>\C[T,T^{-1}] @>>>\frac{\C[T]}{((T-\lambda)^{m+1})}@>>>0\\
@.@V(T-\lambda)VV @VVV @VVV\\
0@>>>\C[T,T^{-1}]@>(T-\lambda)^{m}>>\C[T,T^{-1}] @>>>\frac{\C[T]}{((T-\lambda)^{m})}@>>>0.
\end{CD}\]
As $m\to \infty$, we obtain an inverse system of short exact sequences. It is obvious that $T$ also acts on the inverse system.   Applying $\biG\bullet$, we get a direct system of distinguished triangles. After taking the $\lambda$-eigenspaces of the direct system, we get another direct system of distinguished triangles
\[\cdots\to D_m \to D_{m+1}\to \cdots.\]
By definition, we know $\displaystyle\varinjlim_m (T-\lambda)^m$ kills $\nearbyd{K}{\lambda}$, which means the direct limit of the third entry of $D_m$ is 0; hence 
\[\varinjlim_m D_m=[ \displaystyle\varinjlim_m i^{-1}Rj_*(j^{-1}\d K\otimes f_0^{-1}L^\vee_{(T-\lambda)^m})\stackrel{tr^\infty_{\lambda}}{\arr} \nearbyd{K}{\lambda} \to 0 \stackrel{+1}{\arr}].\]
Since the direct limit functor is exact, we know $\varinjlim D_m$ is a distinguished triangle. Therefore, we see 
\[tr^{\infty}_\lambda : \varinjlim_m i^{-1}Rj_*(j^{-1}\d K \otimes f_0^{-1}L^\vee_{j^{-1} (T-\lambda)^m})[-1]\map \Psi_{f,\lambda}(\d K)[-1]\]
is a quasi-isomorphism (hence an isomorphism as perverse sheaves up to a shift of degree). After replacing $L^\vee_{(T-\lambda)^m}$ by $H_{-\alpha, m}$, the second statement follows by Lemma \ref{dualrr} and isomorphism \eqref{JordanLS}.
\ep 

By the above theorem, in the sense of isomorphisms \eqref{altlambda} and \eqref{altone} we know
\begin{equation} \label{canvarno}
\Psi_{f,\lambda}(\d K)\stackrel{\can}{\simeq}\Phi_{f,\lambda}(\d K)\simeq\varinjlim_m i^{-1}Rj_*(j^{-1}\d K\otimes f_0^{-1}H_{-\alpha,m})
\end{equation}
for $\lambda\neq 1$ and
\begin{equation} \label{canvarone}
\vanishd{K}{1}\simeq \Cone({i^{-1}\d K\map \varinjlim_mi^{-1}Rj_*(j^{-1}\d K\otimes f_0^{-1}H_{-\alpha,m})}).
\end{equation} 
Furthermore, the morphism $\Var$ becomes
\begin{equation}\label{VarJordan}
\Var_\lambda=
\begin{cases}
(0, -2\pi\sqrt{-1}J_{0, \infty}) &\text{ if } \lambda=1,\\
-2\pi\sqrt{-1}J_{0,\infty} &\text{ otherwise,}
\end{cases}
\end{equation}
where $\displaystyle J_{0, \infty}=\varinjlim_m J_{0,m}$. 
 
\section{Nearby and Vanishing Cycles for $\Dmod_X$-modules and Comparisons} 
In this section, we will construct nearby cycles and vanishing cycles for $\Dmod$-modules via $V$-filtrations, and prove comparison theorems. 
\subsection{Koszul complexes and de Rham functor}
First, let us recall the definition of Koszul complexes. Suppose $A$ is an abelian group. Let $\phi_i\in \End_{\Z}A$ be morphisms of $A$ commuting pairwise. Let $K_1=\Cone(A\xmap{\phi_1} A)$ and $K_{i+1}=\Cone(K_i \xmap{\phi_{i+1}}K_i)$; $K(\phi_1,\dots,\phi_n; A)= K_n$ is the Koszul complex of $(A; \phi_1,\dots,\phi_n)$. Then we have the following easy but useful lemma.
\blm\label{Kos}
With $A$ and $\phi_i$ as above, we have:
\be
\item $K(\phi_1,\dots,\phi_n; A)$ is independent of the order of $\phi_i$;
\item If one of $\phi_i$ is an isomorphism, then the Koszul complex is acyclic.
\ee
\elm

Now let $X$ be a complex  manifold of dimension $n$ and $M$ a left $\Dmod_X$-module. Recall that the de Rham complex of $M$ is 
$$
\DR_X(M)=[0\map M\xmap{\nabla}\Omega^1_X\otimes M \map\cdots \map\Omega^n_X\otimes M \map 0]
$$
with $\nabla(m)=\sum_idx_i\otimes\partial_im$ with a local coordinate $(x_1,..., x_n)$.
One can see locally
$$\DR_X(M) = K(\partial_1, \dots, \partial_n; M)[-n].$$ 
Now let $Z=(x_n=0)$, by definition of Koszul complexes we also see 
$$\DR_X(M)\simeq[\DR_Z(M)\xmap{\partial_n}\DR_Z(M)].$$ 

\subsection{V-filtrations on specilizable $\Dmod$-modules} 
Suppose $M$ is a specializable $\Dmod$-module on $X$ along a smooth hypersurface $H$ locally defined by $t=0$. By Corollary \ref{speKM}, the Kashiwara-Malgrange filtration $V^\bullet M$ along $H$ exists. We have the following lemma specializing Lemma \ref{decomp}. 
\blm
Suppose $M$ is specializable along $H$. Then $\gr_V^k M$ has a locally finite decomposition
$$\gr_V^k M=\bigoplus_{\textup{eigenvalues of }t\partial_t} (\gr_V^{k}M)_\alpha,$$
with respect to the $t\partial_t$-action.
\elm

The decompositions in the above lemma give rise to a refinement of $V^\bullet M$ as follows. First, we use the standard order on $\C$; that is for $\alpha_1,\alpha_2\in \C$, $\alpha_1>\alpha_2$ if $\Re(\alpha_1)>\Re(\alpha_2)$ or $\Re(\alpha_1)=\Re(\alpha_2)$ and $\Im(\alpha_1)>\Im(\alpha_2)$. Then for every $\alpha\in \C$, define $V^{\alpha}M$ to be the pre-image of 
$$\displaystyle\bigoplus_{\alpha\le \beta<k+1} (\gr_V^{k}M)_\beta\subseteq \gr_V^k M$$
under the projection $V^k M\arr \gr_V^{k}M$ 
where $k$ is the only interger satisfying $k\le \alpha< k+1$; we define $V^{>\alpha}M$ in the same way, but taking the direct sum over $\alpha< \beta < k+1$.
The resulting filtration $V^\bullet M$ refines the Kashiwara-Malgrange filtration, called the $V$-filtration. Then we can define
\[\gr_V^{\alpha}M=\frac{V^\alpha M}{V^{>\alpha}M}.\]
From construction, we know the $V$-filtration is a $\C$-indexed, locally discrete decreasing filtration satisfying
\be
\item $\bigcup V^\alpha M=M$ and $V^\alpha M$ is coherent over $V^0\Dmod_X$;
\item $tV^\alpha M\subset  V^{\alpha+1}M$, $\partial_tV^\alpha M\subset V^{\alpha-1}M$;
\item $tV^\alpha M = V^{\alpha+1}M$ for $\alpha> -1$;
\item the operator $t\partial_t-\alpha$ acts nilpotently on $\gr_V^{\alpha}M$ for all $\alpha$.
\ee

In particular, the $V$-filtrations exist on holonomic $\Dmod$-modules along every smooth hypersurface. 
\blm
With above notations, we have locally
\be\label{canvarm}
\item $\gr_{V}^{\alpha}M\xmap{t}\gr_{V}^{\alpha+1}M$ is isomorphic if $\alpha\neq -1$; 

\item $\gr_{V}^{\alpha}M\xmap{\partial_t}\gr_{V}^{\alpha-1}M$ is isomorphic if $\alpha\neq 0$.
\ee

\elm
\bp
Both of the statements follow from the nilpotency of the operator $t\partial_t-\alpha$ on $\gr_V^{\alpha}M$.
\ep
From the above lemma, we see $\gr_{V}^{-1}M$ and $\gr_{V}^{0}M$ are distinguished and when $\alpha\notin \Z$, $\gr_{V}^{\alpha}M$ is periodic with period 1 modulo $t$ and $\partial_t$ actions. 

\subsection{Nearby and vanishing cycles for $\Dmod$-modules}
Let $M$ be a coherent $\Dmod_X$-module on a complex manifold $X$. Assume $f$ is a holomorphic function on $X$. We consider the graph embedding 
$$i_f: X\hookrightarrow Y=X\times \C, \text{        } x\mapsto (x, f(x)).$$
Let $t$ be the coordinate on $\C$. Denote $Y^*=Y\setminus X\times \{0\}$ and $X^*=X\setminus f^{-1}(0)$, $i_Y:X\times \{0\} \hookrightarrow Y$ and $i: f^{-1}(0)\hookrightarrow X$ the closed embeddings and $j_Y: Y^*\hookrightarrow Y$ and $j: X^*\hookrightarrow X$ the open embeddings.

\bde
A coherent $\Dmod_X$-module $M$ is specializable along $f$ if $M_f=i_{f,+}M$ is specializable along $X=X\times \{0\}$, where $i_{f,+}$ is the $\Dmod$-module pushforward of $i_f$ (see for instance \cite{Bj} for the definition of pushforward functors for $\Dmod$-module; see also \cite[\S2.4]{BMS} for the explict description of $M_f$). It turns out the $V$-filtration along $X$ on $M_f$ exists. Moreover, $\gr^\alpha_VM_f$ is a coherent $\Dmod_X$-module supported on the divisor $(f=0)$.
\ede
The above definition is compatible with Definition \ref{smspe} because of the following lemma.
\blm
If the hypersurface $H$ defined by $f$ is smooth, then $M$ is specializable along $H$ if and only if $M_f$ is specializable along $X$. Moreover, if $M$ is specializable along $H$, then 
\[ \gr^\alpha_VM_f\simeq  i_+\gr^\alpha_VM\]
for every $\alpha$.
\elm
\bp
Since $f$ is smooth, locally $M_f=M[\partial_t]$. If $M$ is specializable along $H$, one can easily check $V^\alpha M_f=V^\alpha M[\partial_t]$ defines the $V$-filtration of $M_f$. Conversely, if $M_f$ is specializable along $X$, then $V^\alpha M=V^\alpha M_f\cap M$ defines the $V$-filtration of $M$. The second statement follows from Kashiwara's equivalence.
\ep

\bde
If $M$ is specializable along $f$, then the $\alpha$-nearby cycle of $M$ along $f$ is defined to be
$$\Psi_\alpha(M) = \gr^\alpha_VM_f,$$ 
and the  $(\alpha-1)$-vanish cycle 
$$\Phi_{\alpha-1}(M) = \gr^{\alpha-1}_VM_f$$ 
for $0\le\alpha< 1$.
\ede

\subsection{Algebraic localization and specializability}
Assume $N$ is a coherent $\Dmod_Y$-module. Recall that the algebraic localization of $N$ along X is $N(*X)=N\otimes_{\OO}\OO_Y[t^{-1}]$.
\begin{prop}\label{VofM*H}
If $N$ is specializable along $X$, then $N(*X)$ is also specializable along $H$. Moreover, we have
\be
\item $\displaystyle\textup{    } V^{\alpha}N(*X)=\sum_{i\geq 0}V^{\alpha+i}N\otimes_{\OO}\OO_Y(iX);$
\item $\displaystyle\textup{    } \gr_{V}^{\alpha}N(*X)=\sum_{i\geq 0}\gr_{V}^{\alpha+i}N\otimes_{\OO}\OO_Y(iX)$.
\ee
\end{prop}

\bp
First, by Lemma \ref{canvarm} we know 
$$t: \gr^{\alpha}_VM\map \gr^{\alpha+1}_VM$$
is an isomorphism for $\alpha>1$. This implies $V^{\alpha}N(*X)=\sum_{i\geq 0}t^{-i}V^{\alpha+i}N$ is coherent over $V^{0}\Dmod_Y$ for every $\alpha$. Also, it is obvious that $t\partial_t-\alpha$ acts on $\gr_{V}^{\alpha}N(*X)$ nilpotently. Last, since $\partial_t V^\alpha N(*X)+V^\alpha N(*X)=V^{\alpha-1}N(*X)$ for $\alpha<-1$, we know $N(*X)$ is a coherent $\Dmod_Y$-module.
\ep

\begin{cor}\label{gr0M*H}
Assume $N$ is specializable along $X$. Then 
\[\gr_{V}^{\alpha}N(*X)\xmap{t}\gr_{V}^{\alpha+1}N(*X)\]
is an isomorphism for any $\alpha$. Moreover, $$\gr_{V}^{\alpha}N=\gr_{V}^{\alpha}N(*X)$$  for $\alpha>-1$.
\end{cor}

\subsection{Deligne meromorphic extensions of Jordan blocks}
For a complex number $\alpha$, set $$e^\alpha_{k}=t^\alpha\cdot \frac{\log^kt}{k!}.$$ Then we define
\[K^\alpha_{m}=\bigoplus_{k=0}^{m-1}\OO_\C[t^{-1}]e^\alpha_{k}\]
with a naturally defined connection $\nabla$ by requiring
\[\nabla e^\alpha_{k}=\frac{1}{t}(\alpha\cdot e^\alpha_{k}+e^\alpha_{k-1}).\]
We see that the residue of $\nabla$ is $J_{\alpha,m}$ and multivalued $\nabla$-flat sections of $K^\alpha_{m}$ are the $\C$-span of
\[\{e^{-N^\alpha_{m}\log t}\cdot e^\alpha_k\}_{k=0,\dots,m}.\]
Therefore, the local system given by the multivalued $\nabla$-flat sections is exactly $H_{\alpha, m}$, the local system of the Jordan block $J_{\alpha, m}$ (see Example \ref{lsjb}); equivalently,
\[\DR_\C(K^\alpha_{m})|_{\C^*}=H_{-\alpha, m}.\]

By Theorem \ref{altnearby}, we can understand nearby cycles of perverse sheaves using Jordan blocks. Indicated by Riemann-Hilbert correspondence, it is worth to study  
$$N\otimes_\OO K^\alpha_{m}=\bigoplus_{k=0}^{m} N(*X)e^\alpha_{k}.$$ 
\begin{prop}\label{VofJM}
Assume $N$ is specializable along $X$. Then $N\otimes_\OO K^\alpha_{m}$ is also specializable along $X$. Moreover, for every pair $\alpha, \beta\in \C$,
\[\textup{    } \gr_{V}^{\beta}(N\otimes_\OO K^\alpha_m)=\bigoplus_{i\geq 0}^m\gr_{V}^{\beta-\alpha}(N(*X))e^\alpha_i.\]
In particular,
\[\gr_{V}^{0}(N\otimes_\OO K^{-\alpha}_m)=\bigoplus_{i\geq 0}^m\gr_{V}^{\alpha}(N(*X))e^{-\alpha}_i.\]
\end{prop}

\bp
The proof is essentially the same as that of Proposition \ref{VofM*H} and hence left to interested readers.
\ep

Later, we will need the following theorem about Riemann-Hilbert correspondence of algebraic localizations. See \cite{Bj} for the proof. 
\bt \label{Kamcomp}
Assume $f$ is a holomorphic function on $X$. For a regular holonomic $\Dmod_X$-module $M$, the natural morphism 
\[\DR_Y(M_f(*X))\arr Rj_*\DR_X(j^*M)\]
is a quasi-isomorphism. In particular, we have
\[\DR_Y(M_f\otimes_\OO K^{\alpha}_{m})\simeq Rj_*(\DR_{X^*}(M)\otimes_\C f_0^{-1}H_{-\alpha, m}),\]
where $f_0=f|_{X\setminus f^{-1}(0)}$. 
\et

\blm \label{Vmgr}\cite[Proposition 6.4.2]{Bj}
Assume $N$ is specializable along $X$. Then we have
$$ [\DR_X(\gr^0_VN \xmap{\partial_t} \gr^{-1}_VN)] \simeq i^{-1}_Y\DR_Y(N).$$
In particular, if $M$ is a regular holonomic $\Dmod_X$-module and $f$ is a holomorphic function on $X$, then
$$[\DR_X\big(\gr_{V}^{0}(M_f\otimes_\OO K^{\alpha}_m) \xmap{\partial_t} \gr_{V}^{-1}(M_f\otimes_\OO K^{\alpha}_m)\big)] \simeq i^{-1}Rj_*\big (\DR_{X^*}(M)\otimes_\C f_0^{-1}H_{-\alpha, m}\big)$$
\elm


\subsection{Comparison theorems}
Now we can compare nearby and vanishing cycles between regular holonomic $\Dmod$-modules and perverse sheaves.

First, we need a preliminary result about infinite Jordan blocks. 
\blm \label{inftylm}
Let $W$ be a $\C$-vector space, and let $\varphi$ be a $\C$-linear operator on $W$ such that $\varphi-\alpha$ acts on $W$ nilpotently. Set $W_N=\bigoplus_{i=0}^\infty W\otimes e_i$ and define
\[\varphi_\infty(w\otimes e_i)=(\varphi-\alpha)w\otimes e_i+w\otimes e_{i-1}
\] 
for $w\in W$ (assume $e_{-1}=0$). Then $\varphi$ is surjective and $\ker(\varphi_\infty)\simeq W$. 
\elm 
\bp
Define a map $W\arr\ker(\varphi_\infty)$ by
\[w\mapsto \sum_{i\ge 0}(-1)^{i}(\varphi-\alpha)^iw\otimes e_i.\]
Clearly, this map is an isomorphism. 

Since for every $w\in W$ and $j\ge 0$ 
\[\varphi_\infty(\sum_{i\ge j}(-1)^{i-j}(\varphi-\alpha)^{i-j}m\otimes e_i)=w\otimes e_j,\]
surjectivity also follows.
\ep

Now we fix a regular holonomic $\Dmod_X$-module $M$ and a holomorphic function $f$ on $X$. By Proposition \ref{VofJM}, we have
\begin{equation}\label{inftygr}
\varinjlim_m \gr^0_V(M_f\otimes_\OO K^{-\alpha}_m)=\bigoplus_{i=0}^\infty \gr_{V}^{\alpha}(M_f(*X))\otimes e^{-\alpha}_i.
\end{equation}
The nilpotent part of the residue of $K_m^{-\alpha}$ induces an endmorphism $J_{0,\infty}$ of $\varinjlim_m \gr^0_V(M_f\otimes_\OO K^{-\alpha}_m)$ by
\[J_{0,\infty}(m\otimes e_i^{-\alpha})=m\otimes e_{i-1}^{-\alpha}.\]

By equality \eqref{inftygr}, Lemma \ref{inftylm} and Corollary \ref{gr0M*H}, we obtain the following corollary.
\begin{cor} \label{limcomp}
If $M$ is specializable along $f$, then there exists a quasi-isomorphism
\[\gr_V^\alpha(M_f(*X))\xmap\simeq \varinjlim_m[\gr^0_V(M_f\otimes_\OO K^{-\alpha}_m)  \xmap{t\partial_t} \gr^0_V(M_f\otimes_\OO K^{-\alpha}_m)].\]
\end{cor}

Now, we establish everything we need to prove the comparison theorem.
\bt[Kashiwara, Malgrange, Saito] \label{fcompnv}
Assume $M$ is a regular holonomic $\Dmod_X$-module, and $f$ is a holomorphic function on $X$. We have quasi-isomorphisms
\[\DR_X(\Psi_\alpha (M)) \simeq \Psi_{f,\lambda}\big(\DR_X(M)\big)[-1]\]
for $0\le\alpha<1$, and 
\[\DR_X(\Phi_\alpha (M)) \simeq \Phi_{f,\lambda}\big(\DR_X(M)\big)[-1]\]
for $-1\le\alpha<0$, where $\lambda=e^{-2\pi\sqrt{-1}\alpha}$,
and the operator $T$ coming from the monodromy action corresponds to
$$T\simeq \DR_X(e^{-2\pi\sqrt{-1}t\partial_t})=\DR_X(\lambda\cdot e^{-2\pi\sqrt{-1}(t\partial_t-\alpha)})$$ under these isomorphisms (since $t\partial_t-\alpha$ is nilpotent). Moreover, we have an isomorphism of quivers
\[\begin{tikzcd}[every arrow/.append style={shift left}]
 \DR_X\big( \Psi_0(M) \arrow{r}{\partial_t} & \Phi_{-1}(M)\arrow{l}{t}\big)\simeq\big[\Psi_{f, 1}\big(\DR_X(M)\big) \arrow{r}{\can}  &\Phi_{f, 1}\big(\DR_X(M)\big) \arrow{l}{\frac{\Var}{2\pi\sqrt{-1}}}\big][-1]
 \end{tikzcd}\]
for $0\le\alpha<1$.
\et
\bp
To simplify notation, set the complex $$\mathcal A^\alpha_m=[\gr_{V}^{0}(M_f\otimes_\OO K^{-\alpha}_m)  \xmap{t\partial_t} \gr_{V}^{0}(M_f\otimes_\OO K^{-\alpha}_m)]$$
in degrees $-1$ and $0$.

By Theorem \ref{Kamcomp} and Lemma \ref{Vmgr}, we see 
\[\DR_X \mathcal A^\alpha_m\simeq i^{-1}Rj_*(\DR_{X^*}(M)\otimes_\C f_0^{-1}H_{-\alpha, m}\big).\]
Since the de Rham functor and the direct limit functor commute, we have 
\begin{equation}\label{nbinfty}
\DR_X \varinjlim_m \mathcal A^\alpha_m\simeq \Psi_{f,\lambda}\big(\DR_X(M)\big),
\end{equation}
thanks to Theorem \ref{altnearby}. Since $\gr_V^\alpha(M(*X))=\gr_V^\alpha(M)$ for $0<\alpha<1$, we see 
\begin{equation}\label{nbcomp}
\DR_X\big(\gr_V^\alpha(M)\big)[1]\simeq \Psi_{f,\lambda}\big(\DR_X(M)\big),
\end{equation}
by the quasi-isomorphism \eqref{nbinfty} and Corollary \ref{limcomp}.

For $0<\alpha<1$, we have $$\can\circ\var=\log T_u=-2\pi\sqrt{-1}\cdot J_{0,\infty}$$ on $\Psi_{f,\lambda}\big(\DR_X(M)\big)$ by \eqref{VarJordan}. On the other hand, one can check that the operator $J_{0,\infty}$ on $\ker(\varinjlim\mathcal A^\alpha_m)$ corresponds to $-(t\partial_t-\alpha)$ on $\gr^\alpha_V(M_f)$under the quasi-isomorphism in Corollary \ref{limcomp} . Therefore, we conclude
\[\DR_X(t\partial_t-\alpha)\simeq \frac{\log T_u}{2\pi\sqrt{-1}}\]
under the quasi-isomorphism \eqref{nbcomp}. Since $t\partial_t-\alpha$ is locally nilpotent on $gr_V^\alpha(M)$ and the monodromy action $T$ on $\Psi_{t,\lambda}\big(\DR_X(M)\big)$ has only one eigenvalue $\lambda$, we obtain
\[\DR_X(e^{2\pi\sqrt{-1}t\partial_t})\simeq T^{-1}.\]

Since $\Psi_{f,\lambda}\big(\DR_X(M)\big)$ and $\Phi_{f,\lambda}\big(\DR_X(M)\big)$ are canonically the same (see \eqref{canvarno}) and $\partial_t:\gr_V^\alpha(M_f)\map\gr_V^{\alpha-1}(M_f)$ is isomorphic (see Lemma \ref{canvarm}) for $0<\alpha<1$, we also obtain the corresponding statements for vanishing cycles.

Now we are dealing with the last statement. The natural morphism $M_f\arr M_f\otimes_\OO K^{0}_m$ induces a triangle 
\[[\gr^0_VM_f \xmap{\partial_t} \gr^{-1}_VM_f]\stackrel{\iota}{\arr} \varinjlim_m \mathcal A^0_m \arr \Cone(\iota)\xmap{+1}.\]
We recover the triangle
\[i^{-1}\d{K}\arr \nearbyd{K}{1} \xmap{\can} \vanishd{K}{1} \xmap{+1}\]
by applying the de Rham functor on the above triangle, where $\d K=\DR_X(M)$, thanks to Theorem \ref{Kamcomp} and Lemma \ref{Vmgr}. Therefore, we have 
\[\DR_X(\Psi_0(M) \xmap{\partial_t} \Phi_{-1}(M))[1]\stackrel{\simeq}{\arr}[\nearbyd{K}{1} \stackrel{\can}{\arr} \vanishd{K}{1}].\]

By Corollary \ref{limcomp} we know the morphism 
\[\gr^0_V(M_f)\arr\varinjlim_m \gr^0_V(M_f\otimes_\OO K^{0}_m) \]
defined by 
\[m\mapsto \sum_{i\ge 0}(-1)^{i}(t\partial_t)^im\otimes e_i\]
induces the quasi-isomorphism 
\begin{equation}\label{nearbyoneiso}
\gr^0_V(M_f)[1]\xmap\simeq\varinjlim \mathcal A^0_m.
\end{equation} Similarly, the morphism
\[\gr^{-1}_V(M_f)\arr \gr^0_V(M_f)\oplus\varinjlim_m \gr^0_V(M_f\otimes_\OO K^{0}_m)\]
defined by 
\[m\mapsto (m, \sum_{i\ge 1}(-1)^{i}(t\partial_t)^itm\otimes e_i)\]
induces the quasi-isomorphism 
\begin{equation}\label{vanoneiso}
\gr^{-1}_V(M_f)[1]\stackrel{\simeq}{\arr}\Cone{\iota}.
\end{equation}

From endmorphism $J_{0,\infty}: \varinjlim \mathcal A^0_m\arr \varinjlim \mathcal A^0_m$, we get a morphism
\[(0, J_{0,\infty)}: \Cone(\iota)\arr  \varinjlim_m \mathcal A^0_m.\]
One can check
\[[\Phi_{-1}(M) \xmap{t} \Psi_{0}(M)][1]\stackrel{\simeq}{\arr}[\Cone(\iota)\xmap{(0, J_{0,\infty})}  \varinjlim_m \mathcal A^0_m]\]
under the quasi-isomorphisms \eqref{nearbyoneiso} and \eqref{vanoneiso}. Therefore, we have
\[\DR_X(\Phi_{-1}(M) \stackrel{t}{\arr} \Psi_{0}(M))[1]\stackrel{\simeq}{\arr}[\Phi_{f, 1}\big(\DR_X(M)\stackrel{\frac{\Var}{2\pi\sqrt{-1}}}{\arr}\Psi_{f, 1}\big(\DR_X(M)],\]
thanks to \eqref{VarJordan}.
\ep

\end{document}